\documentclass{gtmon_a}
\pdfoutput=1

\proceedingstitle{The Zieschang Gedenkschrift}
\conferencestart{5 September 2007}
\conferenceend{8 September 2007}
\conferencename{Conference in honour of Heiner Zieschang}
\conferencelocation{Toulouse, France}

\editor{Michel Boileau}
\givenname{Michel}
\surname{Boileau}

\editor{Martin Scharlemann}
\givenname{Martin}
\surname{Scharlemann}

\editor{Richard Weidmann}
\givenname{Richard}
\surname{Weidmann}

\title[The rank of the fundamental group of certain hyperbolic 
3--manifolds]{The rank of the fundamental group of certain\\
hyperbolic 3--manifolds fibering over the circle}

\author{Juan Souto}
\givenname{Juan}
\surname{Souto}
\address{Department of Mathematics\\University of Chicago\\\newline
5734 S University Avenue\\Chicago, Ill 60637\\USA}
\email{juan@math.uchicago.edu}
\urladdr{www.math.uchicago.edu/~juan}

\volumenumber{14}
\issuenumber{}
\publicationyear{2008}
\papernumber{22}
\startpage{505}
\endpage{518}

\doi{}
\MR{}
\Zbl{}

\arxivreference{math/0503163}


\received{11 April 2006}
\revised{26 December 2006}
\accepted{2 January 2007}
\published{29 April 2008}
\publishedonline{29 April 2008}
\proposed{}
\seconded{}
\corresponding{}
\version{}

\makeatletter
\def\cnewtheorem#1[#2]#3{\newtheorem{#1}{#3}[section]
\expandafter\let\csname c@#1\endcsname\c@sat}
\makeatother

\AtBeginDocument{\let\bar\wbar\let\tilde\wtilde\let\hat\what}

\newtheorem{sat}{Theorem}[section]
\cnewtheorem{lem}[sat]{Lemma}
\cnewtheorem{kor}[sat]{Corollary}
\cnewtheorem{prop}[sat]{Proposition}
\newtheorem*{sat*}{Theorem}
\newtheorem*{kor*}{Corollary}

\theoremstyle{definition}

\cnewtheorem{bei}[sat]{Example}
\newtheorem*{defi*}{Definition}
\newtheorem*{bei*}{Example}
\newtheorem*{rmk*}{Remark}
\makeautorefname{sat}{Theorem}

\numberwithin{equation}{section}

\theoremstyle{plain}
\newtheorem*{namedtheorem}{\theoremname}
\newcommand{\theoremname}{testing}
\newenvironment{named}[1]{\renewcommand{\theoremname}{#1}
\begin{namedtheorem}}{\end{namedtheorem}}

\theoremstyle{remark}

\newtheorem{Claim}{Claim}

\newcommand{\BC}{\mathbb C}		\newcommand{\BH}{\mathbb H}
\newcommand{\BR}{\mathbb R}		
		
\newcommand{\BS}{\mathbb S}		\newcommand{\BZ}{\mathbb Z}
\newcommand{\BF}{\mathbb F}		

\newcommand{\CN}{\mathcal N}
\newcommand{\CS}{\mathcal S}		
\newcommand{\bnd}{\partial}

\DeclareMathOperator{\PSL}{PSL}		
\DeclareMathOperator{\Isom}{Isom}	
\DeclareMathOperator{\Map}{Map}
\DeclareMathOperator{\inj}{inj}
\DeclareMathOperator{\diam}{diam}
\DeclareMathOperator{\rank}{rank}
\DeclareMathOperator{\Axis}{Axis}

\begin{document}

\begin{asciiabstract}
We determine the rank of the fundamental group of those hyperbolic 
3-manifolds fibering over the circle whose monodromy is a sufficiently 
high power of a pseudo-Anosov map. Moreover, we show that any two 
generating sets with minimal cardinality are Nielsen equivalent.
\end{asciiabstract}

\begin{abstract}
We determine the rank of the fundamental group of those hyperbolic 
3--manifolds fibering over the circle whose monodromy is a sufficiently 
high power of a pseudo-Anosov map. Moreover, we show that any two 
generating sets with minimal cardinality are Nielsen equivalent.
\end{abstract}

\maketitle

\section{Introduction}

Probably the most basic invariant of a finitely generated group is its
{\em rank}, ie the minimal number of elements needed to generate
it. In general the rank of a group is not computable. For instance,
there are examples, due to Baumslag, Miller and Short
\cite{Baumslag-Miller-Short}, of hyperbolics groups showing that there
is no uniform algorithm solving the rank problem. Everything changes
in the setting of 3--manifold groups and recently Kapovich and Weidmann
\cite{Kapovich-Weidmann-rank} gave an algorithm determining
$\rank(\pi_1(M))$ when $M$ is a 3--manifold with hyperbolic
fundamental group. However, it is not possible to give a priori bounds
on the complexity of this algorithm and hence it seems difficult to
use it to obtain precise results in concrete situations. The goal of
this note is to determine the rank of the fundamental group of a
particularly nice class of 3--manifolds.

Let $\Sigma_g$ be the closed (orientable) surface of genus 
$g\ge 2$, $F\co \Sigma_g\to\Sigma_g$ a mapping class and 
$$M(F)=\Sigma_g\times[0,1]/(x,1)\simeq(F(x),0)$$
the corresponding mapping torus. By construction, $\pi_1(M(F))$ is 
a HNN-extension of $\pi_1(\Sigma_g)$ and hence, considering generating sets of $\pi_1(\Sigma_g)$ with 
$$\rank(\pi_1(\Sigma_g))=2g$$
elements and adding a further element corresponding to the extension we obtain generating sets of $\pi_1(M(F))$ with $2g+1$ elements. We will say that the so-obtained generating sets are {\em standard}. In this note we prove:
\begin{sat}\label{main}
Let $\Sigma_g$ be the closed surface of genus $g\ge 2$, $F\in\Map(\Sigma_g)$ a pseudo-Anosov mapping class and $M(F^n)$ the mapping torus of $F^n$. There is $n_F$ such that for all $n\ge n_F$
$$\rank(\pi_1(M(F^n)))=2g+1.$$
Moreover for any such $n$ any generating set of $\pi_1(M(F^n))$ with minimal cardinality is Nielsen equivalent to a standard generating set.
\end{sat}

Recall that two (ordered) generating sets $\CS=(g_1,\dots,g_r)$ and $\CS'=(g_1',\dots,g_r')$ are {\em Nielsen equivalent} if they belong to the same class of the equivalence relation generated by the following three moves:
$$\begin{array}{ll} 
\hbox{Inversion of}\ g_i&\hspace{1pt}
\left\{\begin{array}{ll}g_i'=g_i^{-1} & \\ g_k'=g_k & k\neq
i\end{array}\right. \\ 
\hbox{Permutation of}\ g_i\ \hbox{and}\ g_j\ \hbox{with}\ i\neq j&
\left\{\begin{array}{ll}g_i'=g_j & \\ g_j'=g_i & \\ g_k'=g_k & k\neq
i,j\end{array}\right. \\ 
\hbox{Twist of}\ g_i\ \hbox{by}\ g_j\ \hbox{with}\ i\neq j &\hspace{1pt}
\left\{\begin{array}{ll}g_i'=g_ig_j & \\ g_k'=g_k & k\neq
i\end{array}\right.
\end{array}$$
It is due to Zieschang \cite{Zieschang} that any two generating sets of $\pi_1(\Sigma_g)$ with cardinality $2g$ are Nielsen equivalent. This implies that any two standard generating sets of a mapping torus $M(F)$ are also Nielsen equivalent. We deduce:

\begin{kor}
Let $\Sigma_g$ be the closed surface of genus $g\ge 2$,
$F\in\Map(\Sigma_g)$ a pseudo-Anosov mapping class and $M(F^n)$ the
mapping torus of $F^n$. There is $n_F$ such that any two minimal
generating sets of $M(F^n)$ are Nielsen equivalent for all $n\ge n_F.\qed$
\end{kor}

In \fullref{sec:nielsen} we recall the relation between 
Nielsen equivalence classes of generating sets of the fundamental 
group of a manifold $M$ and free homotopy classes of graphs in $M$. 
Choosing such a graph with minimal length we obtain a link between 
the algebraic problem on the rank of $\pi_1(M)$ and the geometry 
of the manifold. 
In \fullref{sec:meat} we prove \fullref{meat} which is 
essentially a generalization of the fact that paths in hyperbolic 
space $\BH^3$ which consist of large geodesic segments meeting at large angles are quasi-geodesic. Hyperbolic geometry comes into the picture through a theorem of Thurston who proved that the mapping torus $M(F)$ of a pseudo-Anosov mapping class admits a metric of constant negative curvature; equivalently, there is a discrete torsion-free subgroup $\Gamma\subset\PSL_2\BC=\Isom_+(\BH^3)$ with $M(F)$ homeomorphic to $\BH^3/\Gamma$. The geometry of the manifolds $M(F^n)$ is well understood and in \fullref{sec:manifolds} we review very 
briefly some facts needed in \fullref{sec:proofmain} 
to prove \fullref{main}.

The method of proof of \fullref{main} is suggested by the proof of a
result of White \cite{White} who proved that the rank of the
fundamental group of a hyperbolic 3--manifold yields an upper bound
for the injectivity radius. Similar ideas appear also in the work of
Delzant \cite{Delzant} on subgroups of hyperbolic groups with two
generators, in the proof of a recent result of Ian Agol relating rank
and Heegaard genus of some 3--manifolds and in the work of
Kapovich and Weidmann \cite{Kapovich-Weidmann-rank}. It should be said
that in fact most arguments here are found in some form in the papers
of Kapovich and Weidmann and that the main result of this note cannot
come as a surprise to these authors. It should also be mentioned that
a more general result in the spirit of \fullref{main}, but in the
setting of Heegaard splittings, is due to Bachmann and Schleimer
\cite{Bachmann-Schleimer}.

Recently Ian Biringer has obtained, using methods similar to those in this paper, the following extension of \fullref{main}:

\begin{sat*}[Biringer]
For every $\epsilon$ positive, the following holds for all but finitely many examples: If $M$ is a hyperbolic 3--manifold fibering over $\BS^1$ with fiber $\Sigma_g$ and with $\inj(M)\ge\epsilon$ then $\rank(\pi_1(M))=2g+1$ and any two generating sets of $\pi_1(M)$ are Nielsen equivalent.
\end{sat*}

Other related results can be found in Namazi and Souto \cite{Hossein}
and Souto \cite{3rank}.

\medskip
{\bf Acknowledgements}\qua I would like to thank to Ian Agol, Michel
Boileau, Yo'av Moriah and Richard Weidmann for many very helpful and
motivating conversations. I also thank Ian Biringer and the referee
for useful comments enhancing the exposition.
This paper was written while the author was a member of the
Laboratoire de math\'ematiques Emile Picard at the Universit\'e Paul
Sabatier.

\section{Nielsen equivalence of generating sets and carrier graphs}\label{sec:nielsen}

Let $M$ be a hyperbolic 3--manifold.

\begin{defi*} 
A map $f\co X\to M$ of a connected graph $X$ into $M$ is a 
{\em carrier graph} if the homomorphism $f_*\co \pi_1(X)\to\pi_1(M)$ 
is surjective. Two carrier graphs $f\co X\to M$ and $g\co Y\to M$ 
are {\em equivalent} if there is a homotopy equivalence $h\co X\to Y$ 
such that $f$ and $g\circ h$ are free homotopic. 
\end{defi*}

To every generating set $\CS=(g_1,\dots,g_r)$ of $\pi_1(M)$ one 
can associate an equivalence class of carrier graphs as follows: 
Let $\BF_\CS$ be the free non-abelian group generated by the 
set $\CS$, $\phi_\CS\co \BF_\CS\to\pi_1(M)$ the homomorphism 
given by mapping the free bases $\CS\subset\BF_\CS$ to the 
generating set $\CS\subset\pi_1(M)$ and $X_\CS$ a graph with 
$\pi_1(X_\CS)=\BF_\CS$. The homomorphism $\phi_\CS\co \BF_\CS\to\pi_1(M)$ 
determines a free homotopy class of maps $f_\CS\co X_\CS\to M$, ie  
a carrier graph, and any two carrier graphs obtained in this way 
are equivalent. The so determined equivalence class is said to be 
the {\em equivalence class of carrier graphs associated to $\CS$}.

\begin{lem}\label{Nielsen}
Let $\CS$ and $\CS'$ be finite generating sets of $\pi_1(M)$ with the same cardinality. Then the following are equivalent:
\begin{enumerate}
\item $\CS$ and $\CS'$ are Nielsen equivalent.
\item There is a free basis $\bar\CS$ of $\BF_{\CS'}$ with $\CS=\phi_{\CS'}(\bar\CS)$.
\item There is an isomorphism $\psi\co \BF_\CS\to\BF_{\CS'}$ with $\phi_\CS=\phi_{\CS'}\circ\psi$.
\item $\CS$ and $\CS'$ have the same associated equivalence classes of carrier graphs.
\end{enumerate}
\end{lem}
\begin{proof}[Sketch of the proof]
The implications (1) $\Rightarrow$ (2) $\Leftrightarrow$ (3)
$\Leftrightarrow$ (4) are almost tautological. The implication (2)
$\Rightarrow$ (1) follows from a Theorem of Nielsen, who proved that
any two free basis of a free group are Nielsen equivalent (see for
example Collins et al \cite{CGKZ}).
\end{proof}

The natural bijection given by \fullref{Nielsen} between the set of Nielsen equivalence classes of generating sets of $\pi_1(M)$ and the set of equivalence classes of carrier graphs $f\co X\to M$ plays a central role in the proof of \fullref{main}.

\begin{center}
\parbox{11.5cm}{
{\bf Convention}\qua From now on we will only consider generating sets of minimal cardinality. Equivalently, we consider only carrier graphs $f\co X\to M$ with $\rank(\pi_1(X))=\rank(\pi_1(M))$.}
\end{center}

Given a carrier graph $f\co X\to M$ and a path $I$ in $X$ we say that
its length is the length, with respect to the hyperbolic metric, of
the path $f(I)$ in $M$. Measuring the minimal length of a path joining
two points in $X$ we obtain a semi-distance $d_{f\co X\to M}$ on $X$
and we define the {\em length} $l_{f\co X\to M}(X)$ of the carrier
graph $f\co X\to M$ as the sum of the lengths of the edges of $X$ with
respect to $d_{f\co X\to M}$. The semi-distance $d_{f\co X\to M}$
induced on $X$ is not always a distance since there may be some edges
of length $0$ but minimality of the generating set ensures that by
collapsing these edges we obtain an equivalent carrier graph on which
the induced semi-distance is in fact a distance. Moreover, this
collapsing process does not change the length of the carrier
graph. From now on we will assume without further remark that the
semi-distance $d_{f\co X\to M}$ is in fact a distance.

\begin{defi*}
A carrier graph $f\co X\to M$ has {\em minimal length} if 
$$l_{f\co X\to M}(X)\le l_{f'\co X'\to M}(X')$$
for every equivalent carrier graph $f'\co X'\to M$. 
\end{defi*}

If $M$ is closed then it follows from the Arzela--Ascoli Theorem that
every equivalence class of carrier graphs contains a carrier graph
with minimal length:

\begin{lem}\label{minimal}
If $M$ is a closed hyperbolic 3--manifold, then every equivalence class of carrier graphs contains a carrier graph with minimal length. Moreover, every such minimal length carrier graph is trivalent, hence it has $3(\rank(\pi_1(M))-1)$ edges, the image in $M$ of its edges are geodesic segments, the angle between any two adjacent edges is $\frac{2\pi}3$ and every simple closed path in $X$ represents a non-trivial element in $\pi_1(M). \qed$
\end{lem}

See White \cite[Section 2]{White} for a proof of \fullref{minimal}.

\section{Quasi-convex subgraphs}\label{sec:meat}
Recall that a map $\phi\co X_1\to X_2$ between two metric spaces is 
a $(L,A)$--quasi-isometric embedding if 
$$\frac 1Ld_{X_1}(x,y)-A\le d_{X_2}(\phi(x),\phi(y))\le
Ld_{X_1}(x,y)+A$$ for all $x,y\in X_1$. A $(L,A)$--quasi-isometric
embedding $\phi\co \BR\to X$ is said to be a quasi-geodesic. Observe
that a $(L,0)$--quasi-isometric embedding is nothing more than a
$L$--bi-Lipschitz embedding. Before going further, we state here and
for further reference the following well-known fact:

\begin{lem}\label{constant}
There are constants $l_0,A>0$ such that for all $L\ge l_0$ the following holds:
\begin{itemize}
\item Every path in hyperbolic space $\BH^3$ which consists of geodesic 
segments of at least length $L$ and such that all the angles are at 
least $\frac\pi 4$ is a $A$--bi-Lipschitz embedding. 
\item If $K\subset\BH^3$ is convex then every geodesic ray 
$\gamma\co [0,\infty)\to\BH^3$ with $\gamma(0)\in K$ meets the 
boundary $\bnd\CN_L(K)$ of the neighborhood $\CN_L(K)$ of 
radius $L$ around $K$ with at least angle $\frac\pi 4$.\qed
\end{itemize}
\end{lem}

It is surprising that the author didn't find any reference in the literature for the second claim of \fullref{constant}. Here is a proof. Choose $l_0$ and $A$ as in the first claim of the lemma. Up to increasing $l_0$ once we may also assume that the image of every $A$--bi-Lipschitz embedding $\phi\co [0,T]\to\BH^3$ is within at most distance $\frac 12l_0$ of the geodesic segment joining $\phi(0)$ and $\phi(T)$. Given a convex set $K\subset\BH^3$, $L\ge l_0$ and $\gamma$ a ray as in the lemma which exists $\CN_L(K)$ then let $t_0$ be be the unique time with $\gamma(t_0)\in\bnd\CN_L(K)$ and let $p\in K$ be the point closest to $\gamma(t_0)$. If the angle between $\gamma$ and $\bnd\CN_L(K)$ is 
less than $\frac\pi 2$, then the curve obtained by juxtaposition of 
$\gamma[0,t_0]$ and the geodesic segment $[\gamma(t_0),p]$ consists of two geodesic segments of at least length $l_0$ and with a corner with angle at 
least $\frac\pi 4$. In particular, by the first part of the lemma, 
it is an $A$--bi-Lipschitz embedding and hence by the choice of $l_0$ its image is within of the geodesic segment $[\gamma(0),p]$. However, by convexity of $K$ we have that the latter segment is contained in $K$; a contradiction. 

If $f\co X\to M$ is a carrier graph in a hyperbolic 3--manifold 
$M$ we denote by $\tilde f\co \tilde X\to\BH^3$ the lift of 
$f$ to a map between the universal covers of $X$ and $M$. We will 
be mainly interested in manifolds whose fundamental group is not free; in this case, the map $\tilde f$ cannot be an embedding. However, subgraphs of $X$ may well be quasi-isometrically embedded. 

\begin{defi*}
A connected subgraph $Y\subset X$ of a carrier graph $f\co X\to M$ is 
{\em $A$--quasi-convex} for some $A>0$ if:
\begin{itemize}
\item The restriction $\tilde f\vert_{\tilde Y}\co \tilde Y\to\BH^3$ of the map $\tilde f$ to the universal cover $\tilde Y$ of $Y$ is an 
$(A,A)$--quasi-isometric embedding.
\item Every point in $\tilde Y$ is at most at distance $A$ from the axis of some element of $\pi_1(Y)$.
\item The translation length of every element $f_*(\gamma)$ in $\BH^3$ is at least $\frac 1A$ for every $\gamma\in\pi_1(Y)$.
\end{itemize}
\end{defi*}

Recall that a discrete subgroup $G$ of $\PSL_2\BC$ is {\em convex--cocompact} 
if there is a convex $G$--invariant subset $C\subset\BH^3$ of hyperbolic 
space with $C/G$ compact. The smallest such convex subset of $\BH^3$ 
is the {\em convex-hull} $CH(G)$ of $G$ and it is well-known that 
$CH(G)$ is the closure of the union of all axis of elements in $G$.

If $Y$ is a graph and $g\co Y\to M$ is a map whose lift 
$\tilde g\co \tilde Y\to\BH^3$ is a quasi-isometric embedding 
then the image $g_*(\pi_1(Y))$ is a free convex--cocompact subgroup. 
Intuitively, considering $A$--quasi-convex graphs amounts to considering 
uniformly convex--co\-compact free subgroups. More precisely, 
if $Y\subset X$ is $A$--quasi-convex and $\gamma\in\pi_1(Y)$ is non-trivial then the image $\tilde f(\Axis(\gamma))$ is an $(A,A)$--quasi-geodesic and hence it is at uniformly bounded distance of the axis $\Axis(f_*(\gamma))$ of $f_*(\gamma)$. In particular, there is a $d$ depending only on $A$ with 
$$\tilde f(\tilde
Y)\subset\CN_d(CH(f_*(\pi_1(Y))))\subset\CN_{2d}(\tilde f(\tilde Y))$$
This fact, together with the last condition in the definition of
$A$--quasi-convex, implies:

\begin{lem}\label{good-convex}
For all $A$ there is $d$ such that for every hyperbolic manifold $M$
and every $A$--quasi-convex subgraph $Y$ of a minimal length carrier
graph $f\co X\to M$ there is a $f_*(\pi_1(Y))$--invariant convex
subset $\bar C(Y)$ with $$\tilde f(\tilde Y)\subset\bar
C(Y)\subset\CN_d(\tilde f(\tilde Y)),$$ and such that $d_\BH^3(x,\gamma
x)\ge l_0$ for all $x\in\bnd\bar C(Y)$ and $\gamma\in
f_*(\pi_1(Y))$. Here $l_0$ is the constant provided by
\fullref{constant}.\qed
\end{lem}

The following result is the main technical point of the proof of \fullref{main}.

\begin{prop}\label{meat}
For all $A,s>0$ there is $L$ such that whenever $M$ is a hyperbolic 3--manifold, $f\co X\to M$ is a minimal length carrier graph with $s$ edges, $Y_1,\dots,Y_k$ are disjoint connected $A$--quasi-convex subgraphs of $X$ then either
\begin{itemize}
\item $\tilde f\co \tilde X\to\BH^3$ is a quasi-isometric embedding and hence $\pi_1(M)$ is free, or
\item the graph $X\setminus\cup_i Y_i$ contains an edge of at most length $L$.
\end{itemize}
\end{prop}

The author suggests to the reader that he or she proves this proposition him or herself. In fact, a proof by picture takes two not particularly complicated drawings and this is clearly much more economic than the proof written below. 

As mentioned by the referee, \fullref{meat} is a particular case of
the main technical result of Kapovich and Weidmann
\cite{Kapovich-Weidmann-freely} and that it can also be derived from
their \cite[Theorem 2.5]{Kapovich-Weidmann-rank}.

\begin{proof}
Let $l_0$ and $d$ be the constants provided by Lemmas \ref{constant}
and \ref{good-convex}. We are going to show that $\tilde f\co \tilde
X\to\BH^3$ is a quasi-isometric embedding whenever every edge in
$X\setminus\cup Y_i$ has at least length $6l_0+4d$. Seeking a
contradiction, assume that this is not the case.  Then there is an
infinite geodesic ray $\gamma\co [0,\infty)\to\tilde X$ whose image
$\tilde f(\gamma)$ is not a quasi-geodesic. If there is some
$t\in(0,\infty)$ such that $\gamma(t,\infty)$ is disjoint from the
union of the preimages of the graphs $Y_i$, then $\tilde
f(\gamma(t,\infty))$ consists of a perhaps short starting segment and
geodesic segments of length at least $6l_0+4d$ meeting with angle
$\frac{2\pi}3$; \fullref{constant} implies that $\tilde
f(\gamma(t,\infty))$, and hence $\tilde f(\gamma)$, is a
quasi-geodesic ray, contradicting our assumption. Similarly, if there
is $t\in(0,\infty)$ such that $\gamma(t,\infty)$ is contained in a
preimage $\tilde Y_i$ of some $Y_i$ then the assumption that $\tilde
f\vert_{\tilde Y_i}$ is a quasi-isometric embedding implies again that
$\tilde f(\gamma(t,\infty))$ is quasi-geodesic, contradicting again
our assumption. This implies that the curve $\gamma$ has to enter and
leave the union of the preimages of $Y_i$ infinitely often.

Let $a_1<b_1<a_2<b_2<a_3<\dots$ be such that $\gamma(a_j,b_j)$ is contained in and $\gamma(b_j,a_{j+1})$ is disjoint from the preimage of $\cup_i Y_i$ for all $j\ge 1$. Let also $Z_j$ be the component of the preimage of $\cup_iY_i$ containing $\gamma(a_j,b_j)$. For every $j$, the path $\gamma(b_j,a_{j+1})$ consists of edges which by assumption have length at least $6l_0+4d$. Let $\gamma(b_j,c_j)$ be the first edge of this path. We claim that most of the length of $\gamma(b_j,c_j)$ is outside of $\CN_{2l_0}(\bar C(\tilde f(Z_j)))$. 
In fact, by \fullref{good-convex}, every point in the boundary 
of $\CN_{2l_0}(\bar C(\tilde f(Z_j)))$ is at most distance $2l_0+d$ 
from $\tilde f(Z_j)$; the assumption that $f\co X\to M$ is a minimal 
length graph implies that $\tilde f(\gamma(b_j,c_j))$ spends at most $2l_0+d$ time within $\CN_{2l_0}(\bar C(\tilde f(Z_j)))$. Let $b_j^+$ be the exit time. Then $\tilde f(\gamma(b_j^+,c_j))$ is a geodesic segment of at least length $4l_0+3d$ which, by \fullref{constant} 
has at least angle $\frac\pi 4$ with the boundary of 
$\CN_{2l_0}(\bar C(\tilde f(Z_j)))$. A similar discussion applies not when exiting but when entering $\CN_{2l_0}(\bar C(\tilde f(Z_{j+1})$; let $a_{j+1}^-$ be the entry time. 

Setting $I_1=\gamma(a_1,b_1^+)$, $J_1=\gamma(b_1^+,a_2^-)$, $I_2=\gamma(a_2^-,b_2^+)$, $J_2=\gamma(b_2^+,a_3^-)$,... we obtain a decomposition of $\gamma(a_1,\infty)$ in segments with the following properties:
\begin{itemize}
\item $\tilde f(I_j)\subset\CN_{2l_0}(\bar C(Z_j))$ and is $A'$--quasi-geodesic for some $A'$ and all $j$
\item For all $j$, $\tilde f(J_j)$ is a path consisting of geodesic segments of at least length $2l_0$, with at least angle $\frac {2\pi}3$ at the vertices, with endpoints in the boundaries of $\CN_{2l_0}(\bar C(Z_j))$ and $\CN_{2l_0}(\bar C(Z_{j+1}))$ and such that the angles with these boundaries at 
the endpoints are at least $\frac\pi 4$.
\end{itemize}
Before going further we observe that for all $j$ we have 
$\tilde f(a_j^-)\neq\tilde f(b_j^+)$ because the homomorphism 
$f_*\vert_{\pi_1(Y_i)}$ is injective for all $i$. Assume now that 
the distance of $\tilde f(a_j^-)$ and $\tilde f(b_j^+)$ is less 
than $l_0$. Then, by \fullref{good-convex} we have that the images 
in $M$ of $\tilde f(a_j^-)$ and $\tilde f(b_j^+)$, and hence the 
images of the segments $\tilde f(a_j^-,a_j)$ and $\tilde f(b_j,b_j^+)$, 
are different. This implies that we can replace equivariantly the 
segment $\tilde f(a_j^-,a_j)$ by the geodesic segment 
$[\tilde f(a_j^-),\tilde f(b_j^+)]$, getting a new carrier graph 
$f'\co X'\to M$ with length
\begin{align*}
l_{f'\co X'\to M'}(X')&\le l_{f\co X\to M}(X)-l(\tilde f(a_j^-,a_j))+l([\tilde f(a_j^-),\tilde f(b_j^+)])\\ 
&\le l_{f\co X\to M}(X)-2l_0+l_0<l_{f\co X\to M}(X)
\end{align*}
This contradicts the minimality of $l_{f\co X\to M}(X)$ and proves
that the distance between the points $\tilde f(a_j^-)$ and $\tilde
f(b_j^+)$ of $I_j$ is less than $l_0$. Let $I_j'$ be the geodesic
segment joining the endpoints of $I_j$; observe that the length of
this homotopy is bounded by some constant $A''$ because $I_j$ is
$A'$--quasi-geodesic for all $j$. Then the path $\gamma$ is properly
homotopic to the path $\gamma'$ obtained as the juxtaposition of the
segments $I_1'\cup J_1\cup I_2'\cup J_2\cup\dots$. This path consists
now of geodesic segments of at least length $l_0$ and meeting with
angles at least $\frac\pi 4$. \fullref{constant} implies that
$\gamma'$ is a quasi-geodesic. Then the same holds for $\gamma$
because the homotopy from $\gamma$ to $\gamma'$ has at most length
$A''$. This yields the desired contradiction.
\end{proof}

\section{Some facts on the geometry of mapping tori}\label{sec:manifolds}
As mentioned in the introduction, the following is the starting point of our considerations:

\begin{sat*}[Thurston \cite{Thurston}]
Let $\Sigma_g$ be the closed surface of genus $g\ge 2$ and $F\in\Map(\Sigma_g)$ a pseudo-Anosov mapping class. Then the mapping torus 
$$M(F)=\Sigma_g\times[0,1]/(x,1)\simeq(F(x),0)$$
admits a hyperbolic metric.
\end{sat*}

The manifold $M(F)$ fibers over the circle with fiber $\Sigma_g$ and monodromy $F$. Let $\pi\co \pi_1(M(F))\to\BZ$ be the homomorphism given by this fibering and observe that $M(F^n)$ is homeomorphic, and hence isometric by Mostow's rigidity theorem, to the cover of $M(F)$ corresponding to the kernel of the composition of $\pi$ and the canonical homomorphism $\BZ\to\BZ/n\BZ$. Let $M'$ be the infinite cover of $M(F)$ corresponding to the kernel of $\pi$; in the sequel we will always consider $M'$ with the unique hyperbolic metric such that the covering $M'\to M(F)$ is Riemannian. Before going further we observe the following fact that we state here for further reference:

\begin{lem}\label{lifting}
For every $D$ there is $n_D$ such that the following holds for all $n\ge n_D$: Every subset $K\subset M(F^n)$ of diameter at most $D$ lifts homeomorphically to $M'.  \qed$
\end{lem}

Many of the arguments used in the present paper rely on properties of finitely generated subgroups of the fundamental group of $M'$.

\begin{prop}\label{covering}
Every proper subgroup $G$ of $\pi_1(M')\simeq\pi_1(\Sigma_g)$ of rank at most $2g$ is free and convex--cocompact.
\end{prop}
\begin{proof}[Sketch of the proof]
The manifold $M'$ is homeomorphic to $\Sigma_g\times\BR$. In particular, every proper subgroup of $\pi_1(M')\simeq\pi_1(\Sigma_g)$ is either free or isomorphic to the fundamental group of a closed surface which covers $\Sigma_g$ with at least degree 2. Any such surface has genus greater than $g$ and hence its fundamental group has rank greater than $2g$. We have proved that the group $G$ is free. A result due to Thurston in this case and to Agol \cite{Agol} and Calegari--Gabai \cite{Calegari-Gabai} in much more generality asserts that $\BH^3/G$ is homeomorphic to the interior of a handlebody. Now, Canary's generalization of Thurston's covering theorem \cite{Canary-covering} implies that $G$ is convex--cocompact.
\end{proof}

\section[Proof of \ref{main}]{Proof of \fullref{main}}\label{sec:proofmain}

As the kind reader may have deduced from the title of this section, we prove here \fullref{main}. But first, as a warm-up, we show the result of White mentioned in the introduction:

\begin{sat*}[White \cite{White}]
For all $r$ there is $R$ such that every closed hyperbolic 3--manifold $M$ with $\rank(\pi_1(M))\le r$ has $\inj(M)\le R$.
\end{sat*}
\begin{proof}
Let $f\co X\to M$ be a minimal length carrier graph in the class of a minimal generating set of $\pi_1(M)$; observe that $X$ has at most $s=3(r-1)$ edges. Denote by $X^{<t}$ the (possibly empty) subgraph of $X$ consisting of the union of all the edges with length less than $t$. Every simple closed circuit in $X^{<t}$ represents a non-trivial element in $\pi_1(M)$ by \fullref{minimal} and has at most length $3t(r-1)$. In particular, it suffices to show that there is $t_r$ depending only on $r$ such that some component $Y$ of $X^{<t_r}$ is not a tree.

Let $l_0$ be the constant provided by \fullref{constant}. Since $M$ is closed we have that $\pi_1(M)$ is not free and in particular $\tilde f\co \tilde X\to\BH^3$ cannot be a quasi-isometric embedding. In particular, $X^{<l_0}$ is not empty by \fullref{minimal} and \fullref{constant}. If every component $Y$ of $X^{<l_0}$ is a tree then $\diam(\tilde Y)=\diam(Y)\le 3(r-1)l_0$ and hence the map 
$$\tilde f\vert_{\tilde Y}\co \tilde Y\to\BH^3$$ is a
$(3(r-1)l_0,3(r-1)l_0)$--quasi-isometric embedding. We obtain from
\fullref{meat} a constant $l_1=l_1(r)$ depending only on $r$ such that
$X^{<l_0}$ is a proper subgraph of $X^{<l_1}$. If again every
connected component of $X^{<l_1}$ is tree then we get $l_2=l_2(r)$
depending only on $r$ such that $X^{<l_1}$ is a proper subgraph of
$X^{<l_2}$. This process can be repeated at most $3(r-1)$ times since
this is the number of edges in $X$; this concludes the proof of
White's Theorem.
\end{proof}

As we see, the proof of White's Theorem yields in fact that every 
generating set $(g_1,\dots,g_r)$ is Nielsen equivalent to a generating 
set $(g_1',\dots,g_r')$ such that the translation length of $g_1'$ is 
uniformly bounded. The idea of the proof of \fullref{main} is to show 
that every generating set of $\pi_1(M(F^n))$ is Nielsen equivalent to a 
generating set such that the translation lengths of all elements but 1 are 
uniformly bounded.

\begin{named}{\fullref{main}}
Let $\Sigma_g$ be the closed surface of genus $g\ge 2$, $F\in\Map(\Sigma_g)$ a pseudo-Anosov mapping class and $M(F^n)$ the mapping torus of $F^n$. There is $n_F$ such that for all $n\ge n_F$
$\rank(\pi_1(M(F^n)))=2g+1$.
Moreover for any such $n$ any generating set of $\pi_1(M(F^n))$ with 
minimal cardinality is Nielsen equivalent to a standard generating set.
\end{named}
\begin{proof}
For all $n$ let $\CS_n$ be a generating set of $\pi_1(M(F^n))$ with minimal cardinality and $f_n\co X_n\to M(F^n)$ a minimal length carrier graph in the equivalence class determined by $\CS_n$. As remarked in the introduction 
$\rank(\pi_1(M(F^n)))\le 2g+1$
and hence $X_n$ has at most $6g$ edges. As in the proof of White's Theorem, we denote by $X_n^{<t}$ the subgraph of $X_n$ consisting of all the edges of $X_n$ of length less than $t$.

\begin{Claim}For every $D$ there are $n_D$ and $A_D$ such that for every subgraph $Y_n$ of $X_n$ of length less than $D$ and such that the image of $\pi_1(Y_n)$ is convex--cocompact one has: $Y_n$ is $A_D$--quasi-convex for all $n\ge n_D$.
\end{Claim}

\begin{proof}[Proof of Claim 1]
To begin with observe that the injectivity radius of the manifold $M(F^n)$ is bounded from below by $\inj(M(F))$ for all $n$. In particular, the last condition in the definition of $A$--quasi-convex is automatically satisfied for every $A$ with $$A^{-1}\le\inj(M(F)).$$
Seeking a contradiction assume that for some $D$ there are sequences
$A_i,n_i\to\infty$ such that for all $i$ there is a subgraph $Y_{n_i}$
of $X_{n_i}$ which has length less than $D$ and fails to be
$A_i$--quasi-convex but such that $(f_{n_i})_*(\pi_1(Y_{n_i}))$ is
convex--cocompact. Composing the map $f_{n_i}\co X_{n_i}\to
M(F^{n_i})$ with the covering $M(F^{n_i})\to M(F)$ we obtain from
the Arzela--Ascoli Theorem that, up to conjugacy in $\pi_1(M(F))$ and
passing to a subsequence, we may assume that
$(f_{n_i})_*(\pi_1(Y_{n_i}))=(f_{n_j})_*(\pi_1(Y_{n_j}))$ are
conjugated for all $i,j$. In particular, the desired contradiction
follows if we show that the map $\pi_1(Y_{n_i})\to
f_*(\pi_1(Y_{n_i}))$ is an isomorphism.

By \fullref{lifting} there is $i_D$ such that for all $i\ge i_D$ the 
graph $Y_{n_i}$ lifts to $M'$. In particular, we obtain from 
\fullref{covering} that $(f_{n_i})_*(\pi_1(Y_{n_i}))$ is a free 
subgroup of $\pi_1(M')$ which has in particular at most the same 
rank as $\pi_1(Y_{n_i})$. Minimality of the generating set ensures that 
$$\rank((f_{n_i})_*(\pi_1(Y)))=\rank(\pi_1(Y_{n_i})).$$
We are done, since every surjective homomorphism between two free 
groups of the same rank is an isomorphism. 
\end{proof}

We use now an argument similar to the one in the proof of White's Theorem to show:

\begin{Claim}
There are $n_1$ and $t$ such that for all $n\ge n_1$ there is a connected component $Y_n$ of $X_n^{<t}$ such that the image of $\pi_1(Y_n)$ into $\pi_1(M(F^n))$ is not convex--cocompact.
\end{Claim}

\begin{proof}[Proof of Claim 2]
As in the proof of White's Theorem we obtain a first constant $t_1$ such that for all $n$ at least one of the components $Y_{n,t_1}^1,\dots,Y_{n,t_1}^{k(n,t_1)}$ of $X_n^{<t_1}$ is not a tree. If for all $n$ the image of the fundamental group of one of these component fails to be convex--cocompact then are done with $t=t_1$. Assume that there is a subsequence $(n_i)_i$ such that the image of $\pi_1(Y_{n_i,t_1}^j)$ is convex--cocompact for all $j$ and $i$. By claim 1 there is a constant $A_1$ such that $Y_{n_i,t_1}^j$ is $A_1$--quasi-convex for all $i,j$. In particular, we obtain from \fullref{meat} a constant $t_2$ such that $X_{n_i}^{<t_1}$ is a proper subgraph of $X_{n_i}^{<t_2}$ for all $i$. If again the image in of the fundamental group of every connected component of $X_{n_i}^{<t_2}$ is convex--cocompact for infinitely many $i$, say for all $i$, then we can repeat the process. The bound on the number of edges of $X_n$ ensures that after at most $6g$ steps we find the desired subgroup.
\end{proof}

We can now conclude the proof of \fullref{main}. Let $S_n$ be a
generating set of $\pi_1(Y_n)$ where $Y_n$ is the connected subgraph
of $X_n$ provided by claim 2, extend it to a generating set $\bar S_n$
of $X_n$ and let $\bar\CS_n$ be the generating set of $\pi_1(M(F^n))$
obtained as the image of $\bar S_n$ under the homomorphism 
$$(f_n)_*\co \pi_1(X_n)\to\pi_1(M(F^n)).$$
By \fullref{Nielsen}, $\bar\CS_n$ is Nielsen equivalent to the minimal
generating set $\CS_n$ we started with. In particular, $\bar\CS_n$ is
minimal as well. The claim of \fullref{main} follows once we prove
that $\bar\CS_n$ is a standard generating set of $\pi_1(M(F^n))$ and
hence has $2g+1$ elements. Observe that since $\bar S_n$ has
$\rank(\pi_1(M(F^n))\le 2g+1$ elements, it suffices to show that the
generating set $S_n$ of $\pi_1(Y_n)$ has $2g$ elements and that its
image under $(f_n)_*$ generates the subgroup $\pi_1(M')$ of
$\pi_1(M(F^n))$ corresponding to the fiber $\Sigma_g$. This is what we
prove next: The graph $Y_n$ is contained in $X_n^{<t}$, where $t$ is
as in claim 2, and therefore it has at most diameter $6gt$. By
\fullref{lifting} there is $n_1$ such that $Y_n$ lifts to $M'$ for all
$n\ge n_1$; in particular, $\pi_1(Y_n)$ does not surject onto
$\pi_1(M(F^n))$ and hence one has
\begin{equation}\label{eq-rank1}
\rank(\pi_1(Y_n))\le\rank(\pi_1(M(F^n)))-1\le 2g.
\end{equation}
On the other hand, since the image of $\pi_1(Y_n)$ into 
$\pi_1(M')\simeq\pi_1(\Sigma_g)$ is not convex--cocompact we 
deduce from \fullref{covering} that $\pi_1(Y_n)$ surjects 
on $\pi_1(M')$; thus
\begin{equation}\label{eq-rank2}
2g=\rank(\pi_1(M'))\le\rank(\pi_1(Y_n)).
\end{equation}
This concludes the proof of \fullref{main}.
\end{proof}

\bibliographystyle{gtart}
\bibliography{link}

\end{document}